\documentclass[10pt,reqno]{amsart}
\usepackage{amsfonts}
\usepackage{float}
\usepackage{amsmath,amssymb}
\usepackage{lscape}
\usepackage[dvipdfmx]{graphicx, color}
\usepackage[all]{xy}

\def\qed{\hfill $\Box$}
\def\proof{\noindent {\sl Proof} :\;  }

\newcommand{\Lcal}{\mathcal{L}}

\newcommand{\E}{\mathcal{E}}

\newcommand{\R}{\mathbb{R}}

\newcommand\rank{\mbox{\rm rank}\,}

\newcommand\Imag{\mbox{\rm Im} }

\def\qed{\hfill $\Box$}
\def\proof{\noindent {\sl Proof} :\;  }

\def\rd{\partial}

\def\bx{\mbox{\boldmath $x$}}
\def\by{\mbox{\boldmath $y$}}

\def\bb{\mbox{\boldmath $b$}}
\def\bc{\mbox{\boldmath $c$}}

\def\bp{\mbox{\boldmath $p$}}
\def\bq{\mbox{\boldmath $q$}}

\def\bzero{\mathbf{0}}

\def\bbx{\mbox{\tiny$\bx$}}
\def\bby{\mbox{\tiny$\by$}}

\def\bbp{\mbox{\tiny$\bp$}}
\def\bbq{\mbox{\tiny$\bq$}}

\newcommand{\til}[1]{\widetilde{#1}}
\numberwithin{equation}{section}

\newtheorem{thm}{\bf Theorem}[section]
\newtheorem{cor}[thm]{\bf Corollary}
\newtheorem{lem}[thm]{\bf Lemma}
\newtheorem{prop}[thm]{\bf Proposition}
\newtheorem{dfn}[thm]{\bf Definition}

\newtheorem{rem}[thm]{\bf Remark}
\newtheorem{exam}[thm]{\bf Example}

\begin{document}
\title[Statistical manifold with degenerate metric]
{Statistical manifold with degenerate metric}
\author[K.~Kayo]{Kaito Kayo}
\address[K.~Kayo]{M2, Graduate School of Information Science and Technology, Hokkaido University,
Sapporo 060-0814, Japan}
\email{kayo.kaito.b9@elms.hokudai.ac.jp}
%

%
%
\keywords{Codazzi structure, Statistical manifold, Para-complex geometry, Dually flat structure}
\dedicatory{}

%
\begin{abstract}
  A {\it statistical manifold} is a pseudo-Riemannian manifold endowed with a {\it Codazzi structure}.
  This structure plays an important role in Information Geometry and its related fields, e.g., a statistical model admits this structure with the Fisher-Rao metric.
  In practical application, however, the metric may be degenerate, and then this geometric structure is not fully adapted.
  In the present paper, for such cases, we introduce the notice of {\em quasi-Codazzi structure} which consists of a possibly degenerate metric (i.e., symmetric $(0,2)$-tensor) and a pair of  coherent tangent bundles with affine connections.
  This is thought of as an affine differential geometry of Lagrange subbundles of para-Hermitian vector bundles and also as a submanifold theory of para-Hermitian space-form.
  As a special case, the quasi-Codazzi structure with flat connections coincides with the quasi-Hessian structure previously studied by Nakajima-Ohmoto.
  The relation among our quasi-Codazzi structure, quasi-Hessian structure and weak contrast functions generalizes the relation among Codazzi structure, dually flat (i.e., Hessian) structure and contrast functions.

\keywords{First keyword \and Second keyword \and More}
\end{abstract}

\maketitle

{\small
\setcounter{tocdepth}{2}
\tableofcontents
}

\newpage

\section{Introduction}
\label{Intro}

A {\em statistical manifold} is a $C^{\infty}$-manifold $M$ endowed with a {\em Codazzi structure} $(h,\nabla,\nabla^{*})$ where $h$ is a pseudo-Riemannian metric and $\nabla$ and $\nabla^{*}$ are connections satisfying certain properties.
The cubic tensor $C:= \nabla h$ is called the {\em Amari-Chentsov cubic tensor} in information geometry \cite{Amari-Nagaoka, Amari}.
In a special case where $\nabla$ and $\nabla^{*}$ are flat, we call it a dually flat or Hessian manifold.

In application of Information Geometry, such as machine learning, non-convex optimization etc., we often encounter situations where the metric $h$ degenerates along some locus $\Sigma$ on $M$.
It is easy to see that for any $C^{\infty}$-connection $\nabla$ the dual $C^{\infty}$-connection $\nabla^{*}$ does never exist along $\Sigma$ (Proposition \ref{prop:the dual connection does not exist if h is degenerate}).
Therefore, no relevant differential geometry approach to the singular setting had been considered so far.
However, a new breakthrough was recently made by Nakajima-Ohmoto \cite{Nakajima-Ohmoto} updating the definition of dually flat structure of Amari \cite{Amari} to fit the singular case, called {\em quasi-Hessian structure}, and showing its high potential for applications.
A key ingredient is the notion of {\em coherent tangent bundle}, which has been introduced by Saji-Umehara-Yamada \cite{SUY} for studying Riemannian geometry of singular wavefronts.

In the present paper, based on the preceding work \cite{Nakajima-Ohmoto}, we propose a novel generalization of the Codazzi structure to the case where metric $h$ can be degenerate.
Key materials are as follows.
Let $M$ be an $n$-dimensional manifold, $E \rightarrow M$ a vector bundle of $\rank 2n$ which has a {\em para-Hermitian structure}, and $\Phi:TM\rightarrow E$ a vector bundle map such that $\Phi(T_{p}M)$ is a Lagrange subspace at every $p\in M$.
As for the involution $I:E\rightarrow E$ of the para-Hermitian structure, let $E^{+}$ and $E^{-}$ be the subbundle of eigenvectors for eigenvalues $+1$ and $-1$, respectively.
Then $E=E^{+}\oplus E^{-}$, and let $\Phi^{\pm}$ denote the projection of $\Phi$ onto $E^{\pm}$ \cite{V.P.P.M}.
Further, we impose $E^{\pm}$ to have affine connections $\nabla^{\pm}$ which satisfy a certain duality.
These data form two coherent tangent bundles on $M$,
\begin{equation*}
  (E^{+}, \Phi^{+}:TM \rightarrow E^{+},\nabla^{+}) \quad \text{and} \quad (E^{-}, \Phi^{-}:TM \rightarrow E^{-},\nabla^{-}).
\end{equation*}
The pseudo-Riemannian metric $\tau$ of $E$ of type $(n,n)$ induces a possibly degenerate symmetric $(0,2)$-tensor $h:= \Phi^{*}\tau$ on $TM$.
We call the pair of $h$ and these data a {\em quasi-Codazzi structure} on $M$ (Definition \ref{def of quasi-Codazzi structure}).
If we are given a Codazzi structure $(h,\nabla,\nabla^{*})$ on $M$, we set $E=TM\oplus T^{*}M$ (Whitney sum),  $\Phi^{+}:=id_{TM}$, $\Phi^{-}$ to be defined by the pairing induced by $h$, and take mutually dual connections by $\nabla$ and $\nabla^{*}$.
Then we naturally obtain a quasi-Codazzi structure (Example \ref{canonical quasi-Codazzi}).
Conversely, a quasi-Codazzi structure with non-degenerate $h$ can be identified with a Codazzi structure (Proposition \ref{prop:quasi-Codazzi structure with non-degenerate => Codazzi structure}).
Emphasized is that our quasi-Codazzi structure is about a refined version of differential geometry of {\em a Lagrange subbundle of a para-Hermitian vector bundle over $M$} and it also leads to a submanifold theory in para-Hermitian space-form.
In fact, a simplest example is any Lagrange submanifold of $T^{*}N$ with some connections on $N$.
In a different context, a nice non-trivial example of quasi-Codazzi structure is found in the work of Saji-Umehara-Yamada \cite{SUY-F2, SUY-F1} on a wavefront in an ambient manifold with {\em constant curvature}.

We will generalize several basic properties of the Codazzi structure to our quasi-Codazzi case (Proposition \ref{prop of quasi-Codazzi str}).
An imposed condition on $\Phi^{\pm}$, the relative torsion-freeness, takes an important role, and especially, a symmetric cubic tensor $C$ is naturally defined.
There is a well-known relation between statistical manifolds and contrast functions in the sense of Eguchi \cite{Eguchi} (and also yokes in Barndorff-Nielsen \cite{Barndorff-Nielsen}).
\begin{table}
  \label{table1}
  \begin{tabular}{|c|c|}
    \hline
    Codazzi structure $(M,h,\nabla,\nabla^{*})$ & ?\\
    \hline
    $(M,h,C)$, & $(M,h,C)$,\\
    where $h$ is a metric and & where $h$ is (0,2)-tensor and\\
    $C$ is a symmetric cubic tensor & $C$ is a symmetric cubic  tensor\\
    \hline
    Dually flat (Hessian) structure & Quasi-Hessian structure\\
    \hline
    contrast function & weak contrast function\\
    \hline
  \end{tabular}
  \caption{}
\end{table}
That is, given a contrast function, we can produce $h$ and $C$ (and hence $\nabla$ and $\nabla^{*}$) in a canonical way, and conversely, for a statistical manifold, we can find a contrast function which recovers $h$ and $C$, but it is not unique.
This relation works even in the weaker setting with possibly degenerate $h$ and {\em weak} contrast functions \cite{Matumoto}.
Our quasi-Codazzi structure $(h,(E,\tau,I),\Phi,\nabla^{+},\nabla^{-})$ fills up the missing blank in the above Table \ref{table1} in a proper sense.
We remark that our quasi-Codazzi structure defines $(M,h,C)$, but the converse may not be true.

Finally, we characterize the quasi-Hessian structure of Nakajima-Ohmoto \cite{Nakajima-Ohmoto} as the quasi-Codazzi structure with flat connections $\nabla^{\pm}$ (Theorem \ref{thm of quasi-Hessian and quasi-Codazzi}).
In fact, the flatness condition implies the integrability of the Lagrange subbundle to construct a Lagrange submanifold in a para-Hermitian space-form.

The rest of the present paper is organized as follows.
In \S2, we give a brief summary on some basics in the theory of statistical manifolds and contrast functions and study several fundamental properties.
In \S3, we introduce the definition of quasi-Codazzi structure precisely.
In \S4, we characterize the quasi-Hessian structure as a special case of quasi-Codazzi structure.

\vskip\baselineskip
\vskip\baselineskip
\section{Statistical manifolds}
\label{sec2}
In this section, we summarize the theory of statistical manifolds.
For more details, see \cite{Amari-Nagaoka, Amari, Matsuzoe}.
Let $M$ be an n-dimensional $C^{\infty}$-manifold with a pseudo-Riemannian metric $h$ and an affine connection $\nabla$ on $TM$.

\vskip\baselineskip
\subsection{Codazzi structure}
\label{sec:2-1}
We define the {\em dual connection} $\nabla^{*}$ of $\nabla$ with respect to $h$ by
\begin{align}
Xh(Y,Z) = h(\nabla_{X}Y, Z) + h(Y,\nabla^{*}_{X}Z) \label{2-a}
\end{align}
for $X,Y,Z \in \Gamma(TM)$.
Since $h$ is non-degenerate, $\nabla^{*}$ is uniquely determined and it holds that $(\nabla^{*})^{*} = \nabla$.
Denote by $R$ the curvature-tensor of $\nabla$, and by $T$ the torsion-tensor of $\nabla$, that is,
\begin{align*}
R(X,Y)Z &:= \nabla_{X} \nabla_{Y}Z - \nabla_{Y}  \nabla_{X}Z - \nabla_{[X,Y]}Z,\\
T(X,Y) &:= \nabla_{X}Y - \nabla_{Y}X - [X,Y].
\end{align*}
When $\nabla$ is curvature-free and torsion-free, we call $\nabla$ a flat connection.
The curvature and the torsion of $\nabla^{*}$ are denoted by $R^{*}$ and $T^{*}$, respectively.
From straightforward calculations, we see
\[
h(R(X,Y)Z,W) = -h(Z,R^{*}(X,Y)W)
\]
for $W,X,Y,Z \in \Gamma(TM)$.
In particular, if $R \equiv 0$, then $R^{*}\equiv 0$, and vice-versa.

\begin{prop}\upshape
  (Matsuzoe \cite[Proposition2.2]{Matsuzoe})
  \label{prop of Codazzi str}
  Consider the following conditions:
  \begin{itemize}
    \item[(i)]
    $\nabla$ is torsion-free (i.e., $T\equiv 0$);
    \item[(ii)]
    $\nabla^{*}$ is torsion-free (i.e., $T^{*}\equiv 0$);
    \item[(iii)]
    $C:=\nabla h$ is totally symmetric;
    \item[(iv)]
    $\nabla^{(0)}:=(\nabla+\nabla^{*})/2$ is the Levi-Civita connection with respect to $h$.
  \end{itemize}
  Then any two conditions imply the rest of them.
\end{prop}

By this proposition, we get a totally symmetric cubic tensor $C$, if affine connections $\nabla$ and $\nabla^{*}$ are mutually torsion-free.
We call $C$  the {\em Amari-Chentsov cubic tensor}.
Conversely, given a totally symmetric cubic tensor $C$, we can define mutually dual affine connections.

\begin{prop}\upshape \cite[Proposition2.3]{Matsuzoe}
  \label{Codazzi str from C}
  Let $C$ be a totally symmetric cubic tensor on $M$, and $\nabla^{(0)}$ the Levi-Civita connection with respect to $h$.
  Set $\nabla$ and $\nabla^{*}$ by
  \begin{align*}
    h(\nabla_{X}Y,Z)&:=h(\nabla^{(0)}_{X}Y,Z) - \frac{1}{2}C(X,Y,Z),\\
    h(\nabla^{*}_{X}Y,Z)&:=h(\nabla^{(0)}_{X}Y,Z)+ \frac{1}{2}C(X,Y,Z).
  \end{align*}
  Then $\nabla$ and $\nabla^{*}$ are mutually dual connections and torsion-free, respectively.
  In addition, it holds that $\nabla h=C$.
\end{prop}

\begin{dfn}\upshape \cite{Amari-Nagaoka, Amari, Matsuzoe}
  If all properties in Proposition $2.1$ hold, we call the triplet $(h,\nabla,\nabla^{*})$ a {\em Codazzi structure} on $M$, and also we call $M$ with this structure a {\em statistical manifold}.
  In particular, if $\nabla$ and $\nabla^{*}$ are flat, then we call $M$ a {\em dually flat} (or {\em Hessian}) {\em manifold}.
\end{dfn}

Given a pseudo-Riemannian metric $h$ and a totally symmetric cubic tensor $C$, we get a statistical manifold $(M,h,\nabla)$ by Proposition \ref{Codazzi str from C}.
So, we may rewrite that it by triplet $(M,h,C)$.

Given a pseudo-Riemannian metric $h$ and a torsion-free connection $\nabla$, it is well known \cite[Theorem 2]{Delanoe} (\cite[Theorem 2.5.1]{Shima}) that $C=\nabla h$ is symmetric, i.e., $(M,h,C)$ is a statistical manifold, if and only if the horizontal space determined by $\nabla$ is Lagrangian in $T(T^{*}M)$ via $\til{h}:TM \ni X \mapsto h(X,-)\in T^{*}M$.
Namely, the Codazzi structure implicitly contains the geometry of a certain Lagrange subbundle.
This observation leads to our definition of quasi-Codazzi structure which will be introduced in the next section.

\begin{exam}\upshape
  Let $(\Omega, \beta)$ be a measurable space.
  A parametric statistical model on $\Omega$ is a set $M$ parametrized by $\zeta = (\zeta^{1}, \cdots \zeta^{n}) \in U \subset \R^{n}$:
  \[
  M= \left\{ p(x,\zeta) \;|\; p(x,\zeta)>0, \;\; \int_{\Omega} p(x,\zeta) \; dx = 1 \right\}.
  \]
  The statistical model $M$ is a manifold with a local coordinate system $\zeta = (\zeta^{1}, \cdots \zeta^{n})$, if it satisfies suitable  conditions (see Chapter $2.1$ in \cite{Amari-Nagaoka}).

  For simplicity, set $l(x,\zeta):= \log p(x,\zeta)$.
  We define a symmetric $(0,2)$-tensor $g$ and a totally symmetric cubic tensor $C$ by
  \begin{align*}
    \textstyle{
    g  \left( \frac{\rd}{\rd \zeta^{i}}, \frac{\rd}{\rd \zeta^{j}} \right)}
    &:= \int_{\Omega} \textstyle{\left( \frac{\rd}{\rd \zeta^{i}} l(x,\zeta) \right)
    \left( \frac{\rd }{\rd \zeta^{j}} l(x,\zeta) \right) p(x,\zeta) \; dx},\\
    C \textstyle{\left( \frac{\rd}{\rd \zeta^{i}}, \frac{\rd}{\rd \zeta^{j}}, \frac{\rd}{\rd \zeta^{k}} \right)} &:= \int_{\Omega} \textstyle{\left( \frac{\rd}{\rd \zeta^{i}} l(x,\zeta) \right) \left( \frac{\rd}{\rd \zeta^{j}} l(x,\zeta) \right) \left( \frac{\rd}{\rd \zeta^{k}} l(x,\zeta) \right) p(x,\zeta) \; dx}.
  \end{align*}
  If $g$ is non-degenerate, $g$ defines a pseudo-Riemannian metric on $M$.
  We call $g$ the {\em Fisher-Rao metric} associated to the statistical model $M$.

  For an arbitrary constant $\alpha \in  \R$, we define a torsion-free connection $\nabla^{\alpha}$ on $M$ by
  \begin{gather*}
    \Gamma^{(\alpha)}_{ij,k} (\zeta):= \int_{\Omega} \textstyle{\left\{ \frac{\rd^{2}}{\rd \zeta^{i} \rd \zeta^{j}}l(x,\zeta) + \frac{1-\alpha}{2} \left( \frac{\rd}{\rd \zeta^{i}}l(x,\zeta) \right) \left( \frac{\rd}{\rd \zeta^{j}}l(x,\zeta) \right)  \right\} \left( \frac{\rd}{\rd \zeta^{k}}l(x,\zeta) \right) p(x,\zeta) \; dx},\\
    g \textstyle{\left(\nabla^{(\alpha)}_{ \frac{\rd}{\rd \zeta^{i}} } \frac{\rd}{\rd \zeta^{j}}, \frac{\rd}{\rd \zeta^{k}}
    \right) := \Gamma^{(\alpha)}_{ij,k}(\zeta)}.
  \end{gather*}
  The affine connection $\nabla^{(\alpha)}$ is called the $\alpha$-{\em connection} on $M$.
  We can check easily that $(\nabla^{(\alpha)})^{*} = \nabla^{(-\alpha)}$ and $\nabla^{(0)}$ is the Levi-Civita connection with respect to $g$.
  In addition, $\nabla^{(\alpha)}g$ and $\nabla^{(-\alpha)}g$ are totally symmetric.
  Then $(g,\nabla^{(\alpha)},\nabla^{(-\alpha)})$ is a Codazzi structure on $M$.
\end{exam}

If $h$ is degenerate, there does not exist a dual connection by the following Proposition.
\begin{prop}\upshape
  \label{prop:the dual connection does not exist if h is degenerate}
  Suppose that $h$ is degenerate at a point $p\in M$ and there exists a $C^{\infty}$-regular curve $c(t)$ on $M$ (i.e., $\frac{dc}{dt} \neq 0$) such that $c(0)=p$ and $h$ is non-degenerate at every point $c(t)$ with $t>0$.
  Then, for any $C^{\infty}$-connections $\nabla$ on $M$, no $C^{\infty}$-connection $\nabla^{*}$ on $M$ satisfies the equality:
  \begin{align*}
      Xh(Y,Z) = h(\nabla_{X}Y,Z) + h(Y,\nabla^{*}_{X}Z)
  \end{align*}
  for $X,Y,Z \in \Gamma(TM)$ at $p$.
\end{prop}

\proof
Suppose such $\nabla^{*}$ exists.
Let $(U;x_{1},\cdots x_{n})$ be an arbitrary neighborhood around $p$ and $ (e_{1},\cdots, e_{n})$ an orthogonal frame on $U$ such that $h(e_{i},e_{j})=\delta_{ij} \lambda_{i}$ with some $\lambda_{i}:M\rightarrow \R$ ($\delta_{ij}$ is the Kronecker symbol).
We set $\omega_{i}^{j}$ (resp. $\theta_{i}^{j}$) are connection forms of $\nabla$ (resp. $\nabla^{*}$) with respect to the frame $(e_{i},\cdots,e_{n})$ which satisfy that:
\[
\nabla_{X}e_{i} = \sum_{j=1}^{n} \omega_{i}^{j}(X)e_{j}, \qquad \nabla^{*}_{X}e_{i}=\sum_{j=1}^{n}\theta_{i}^{j}(X)e_{j}
\]
with $X\in \Gamma(TU)$.
Then
\begin{align*}
  0 &= Xh(e_{i},e_{j}) - h(\nabla_{X}e_{i},e_{j}) - h(e_{i},\nabla^{*}_{X}e_{j})\\
  &=X(\delta_{ij}\lambda_{i}) - \sum_{l=1}^{n} \omega_{i}^{l}(X)h(e_{l},e_{j}) -\sum_{l=1}^{n} \theta_{j}^{l}(X)h(e_{i},e_{l})\\
  &= \delta_{ij}X\lambda_{i} - \omega_{i}^{j}(X)\lambda_{j} - \theta_{j}^{i}(X)\lambda_{i}.
\end{align*}
In particular, we get
\begin{align*}
X\lambda_{i} - \omega_{i}^{i}(X) \lambda_{i} - \theta_{i}^{i}(X)\lambda_{i}=0.
\end{align*}
We assume that $\lambda_{1}(p)=0$ at $p\in U$ and $(X\lambda_{1})(p)=0$ for every $X\in \Gamma(TU)$.
Then
\[
(\omega_{1}^{1}(X) + \theta_{1}^{1}(X))(c(t)) =
  \displaystyle{\frac{(X\lambda_{1}) (c(t))}{\lambda_{1} (c(t))}} \qquad (t > 0)
\]
by the assumption $\lambda_{1} \circ c(t) \neq 0$ ($t>0$).
Put $f(t):= \lambda_{1} \circ c(t)$ and $X = \frac{dc}{dt}$ (i.e., $f(0)=0$, $f^{\prime}(0)=0$ and $Xf(t)=f^{\prime}(t) \neq 0$ $(t>0)$).
Then $| f^{\prime}(t)/f(t)| \rightarrow \infty$ at $t \rightarrow 0+0$ (otherwise, $\log |f(t)|$ is bounded, but $f(0)=0$).
Thus $\omega^{1}_{1} + \theta^{1}_{1}$ does not exist at $p$.
This makes a contradiction.
\qed

\vskip\baselineskip
\subsection{Contrast functions}
\label{sec:2-2}

Let $\rho:U \rightarrow \R$ be a function defined on an open neighborhood $U$ of the diagonal $\Delta_{M} \subset M\times M$.
A function $\rho[X_{1} \cdots X_{k}|Y_{1} \cdots Y_{l}]$ is defined by
\[
(X_{1})_{p} \cdots (X_{k})_{p} (Y_{1})_{q} \cdots (Y_{l})_{q} (\rho(p,q)) |_{p=q=r}
\]
for every $r \in M$ and $X_{i},Y_{j}\in \Gamma(TM)$ $(1 \leq i \leq k, 1 \leq j \leq l)$.
We also define $\rho[X|-](r) = X_{p}\rho(p,q)|_{p=q=r}$ and $\rho[-|X](r) = X_{q}\rho(p,q)|_{p=q=r}$.

We call $\rho$ a {\it contrast function} if it satisfies the following conditions:
\[
\text{(i)}\rho[-|-] = \rho(p,p) = 0; \qquad
\text{(ii)}\rho[X|-] = \rho[-|X]=0;
\]
\begin{center}
  (iii)$h(X,Y):=-\rho[X|Y]$ is a pseudo-Riemannian metric on $M$.
\end{center}
If it satisfies only (i) and (ii), $\rho$ is called a {\it weak contrast function}.

When $\rho$ is a contrast function, we define affine connections by
\[
h(\nabla_{X}Y, Z) := -\rho[XY|Z], \qquad
h(Y, \nabla^{*}_{X}Z) := -\rho[Y|XZ].
\]
Those connections are well-defined since $h$ is non-degenerate.
Furthermore, $\nabla$ and $\nabla^{*}$ are torsion-free and mutually dual affine connections with respect to $h$.
Therefore $(M, h,\nabla)$ is a statistical manifold.

Conversely, given a Codazzi structure, we can find a contrast function which reproduces the structure.

\begin{thm}\upshape\cite{Matumoto}
  \label{weak contrast function}
  Given a symmetric $(0,2)$-tensor $h$, i.e., a possibly degenerate metric, and a totally symmetric cubic tensor $C$, we get a weak contrast function $\rho$ such that, for any $X,Y,Z \in \Gamma(TM)$,
  \begin{align*}
    h(X,Y) &=-\rho[X|Y] \; (=\rho[XY|-] = \rho[-|XY]), \\
    C(X,Y,Z) &=-\rho[Z|XY] + \rho[XY|Z].
  \end{align*}
\end{thm}

\vskip\baselineskip
\vskip\baselineskip
\section{Statistical manifolds with degenerate metric}
\label{sec:3}
In this section, we propose a generalization of statistical manifolds so that it admits the metric to be degenerate and fits the theory of weak contrast functions as indicated in Theorem \ref{weak contrast function}.
Naively, we may have expected that there exists a connection $\nabla^{*}$ on $M$ determined by the equality \eqref{2-a} for a given connection $\nabla$, even if $h$ is degenerate along some locus.
However, it is not true: such $\nabla^{*}$ does never exist as a $C^{\infty}$-connection by Proposition \ref{prop:the dual connection does not exist if h is degenerate}.
As a conclusion, we definitely need a new framework to describe the geometry of Codazzi structure under a singular setting.
We employ an affine geometry version of coherent tangent bundles due to Saji-Umehara-Yamada \cite{SUY}.

\vskip\baselineskip
\subsection{Coherent tangent bundles}
\label{sec:3-1}
Let $M$ be an $n$-dimensional $C^{\infty}$-manifold, $\mathcal{E}$ a vector bundle of $\rank n$ on $M$, $\nabla$ a connection on $\mathcal{E}$ and $\Phi : TM \rightarrow \mathcal{E}$ a bundle map.

\[
\xymatrix{
TM \ar[rd]_{\pi_{TM}} \ar[rr]^{\Phi} & & \E \ar[ld]^{\pi_{\E}} \\
& M &
}
\]

\begin{dfn}\upshape
  \label{def of relative torsion and coherent tangent bundle}
  (cf. \cite{SUY})
  We define a {\em relative torsion} $T^{\nabla}:\Gamma(TM) \times \Gamma(TM) \rightarrow \Gamma(\mathcal{E})$ with respect to $\Phi$ by
  \begin{align*}
    T^{\nabla}(X,Y) := \nabla_{X}\Phi(Y) - \nabla_{Y}\Phi(X)- \Phi([X,Y])
  \end{align*}
  for $X,Y\in \Gamma(TM)$.
  The triplet $(\E,\Phi,\nabla)$ is called a {\em coherent tangent bundle} on $M$ if $\nabla$ is relatively torsion-free, i.e.,  $T^{\nabla}\equiv0$.
\end{dfn}

\vskip\baselineskip
\subsection{Para-complex structure}
\label{sec:3-2}
To begin with, we summarize para-complex geometry; the readers are referred to \cite{V.P.P.M} for the detail.

Let $M$ be a $2n$-dimensional $C^{\infty}$-manifold.
An $(1,1)$-tensor $I:TM \rightarrow TM$ with $I^{2}=id$ is called an {\em almost product structure}.
Then, $I$ has eigenvalues $1$ and $-1$; we denote by $T^{+}M$ (resp. $T^{-}M$) the subbundle of eigenvectors for the eigenvalue $1$ (resp. $-1$).
If $T^{+}M$ and $T^{-}M$ have the same $\rank$, $I$ is called an {\em almost para-complex structure} and $(M, I)$ is called an {\em almost para-complex manifold}.
Moreover, if $I$ is integrable, i.e, $[X,Y]\in \Gamma(T^{+}M)$ for $X,Y \in  \Gamma(T^{+}M)$, then $I$ is called a {\em para-complex structure} and $(M,I)$ is called a {\em para-complex  manifold}.

Let $(M,I)$ be a para-complex manifold and $\tau$ a pseudo-Riemannian metric on $M$.
The triplet $(M,\tau,I)$ is called a {\it para-Hermitian manifold} if it holds that
\begin{align}
  \tau (IX,Y) + \tau (X,IY) = 0 \label{def of para-Hermitian}
\end{align}
for every $X,Y\in \Gamma(TM)$.
Then, it is easy to check $\tau$ has signature $(n,n)$.

We define a skew-symmetric $(0,2)$-tensor field $\omega$ by
\[
\omega(X,Y) := \tau(X,IY).
\]
This is actually a symplectic form on $TM$.
If $\omega$ is closed, i.e., $d\omega = 0$, $(M,\tau, I)$ is called a {\it para-K\"{a}hler manifold}.

\begin{exam} \upshape \label{ex para-Kahler}
  Let $(x_{1}, \cdots x_{n}, y_{1}, \cdots, y_{n})$ be the standard basis on $\R^{2n}$.
  An $(1,1)$-tensor $I$ is defined by
  \[
  \textstyle{I \left( \frac{\partial}{\partial x_{i}} \right) = \frac{\partial}{\partial x_{i}}}, \quad \textstyle{I \left( \frac{\partial}{\partial y_{i}} \right) = -\frac{\partial}{\partial y_{i}}}
  \]
  for $i \in \{ 1, \cdots ,n\}$.
  Then $(\R^{2n},I)$ is a para-complex manifold.
  Furthermore, we consider a pseudo-Riemannian metric $\tau$ such that
  \[
    \textstyle
    \tau \left( \frac{\partial}{\partial x_{i}}, \frac{\partial}{\partial x_{j}} \right) =
    \tau \left( \frac{\partial}{\partial y_{i}}, \frac{\partial}{\partial y_{j}} \right) =0,
    \quad
    \tau \left( \frac{\partial}{\partial x_{i}}, \frac{\partial}{\partial y_{j}} \right) = \delta_{ij}
  \]
  for every $i,j\in \{ 1,\cdots, n \}$ and $\delta_{ij}$ is the Kronecker symbol.
  Obviously, this structure satisfies the equation \eqref{def of para-Hermitian}.
  Therefore $(\R^{2n}, \tau, I)$ is a para-Hermitian manifold.
  In addition, the symplectic form $\omega$ is written just as $\omega = -\Sigma_{i=1}^{n} dx_{i} \wedge dy_{i}$ and satisfies $d\omega = 0$.
  Thus $(\R^{2n},\tau, I)$ is a para-K\"{a}hler manifold.
\end{exam}

Let $(M,I)$ and $(M^{\prime},I^{\prime})$ be $2n$-dimensional para-complex manifolds.
A smooth map $f:M \rightarrow M^{\prime}$ is called a {\em para-holomorphic} if $df \circ I = I^{\prime} \circ df$.


These basic knowledge of para-complex geometry is defined on tangent bundle.
Let us consider para-complex geometry for vector bundles.
The following structure is important in the present paper.

\begin{dfn}\upshape
  \cite[Definition 3]{LS05}.
  Let $M$ be an $n$-dimensional $C^{\infty}$-manifold, $E$ a vector bundle of $\rank 2n$ on $M$, $\tau$ a pseudo-Riemannian metric of the bundle $E$, and $I:E \rightarrow E$ a bundle map.
  We call $(E,\tau,I)$ a para-Hermitian vector bundle on $M$ if the following conditions hold:
  \begin{itemize}
    \item[(i)]
    $I^{2}=id$ and $I\neq id$;
    \item[(ii)]
    Let $E^{+}$ (resp. $E^{-}$) be the subbundle of $E$ consisting of eigenvectors with the eigenbundle $1$ (resp. $-1$).
    Then $E^{+}$ and $E^{-}$ have the same rank and $E=E^{+}\oplus E^{-}$;
    \item[(iii)]
    $\tau(IX,Y)+\tau(X,IY)=0$.
  \end{itemize}
\end{dfn}

\begin{rem}\upshape
  Set $\omega(X,Y):=\tau(X,IY)$ to be a symplectic form on $E$.
  And, by the condition (iii), $\tau|_{E^+}=0$ and $\tau|_{E^{-}}=0$.
\end{rem}

\begin{dfn}\upshape
  Let $(E,\tau,I)$ and $(E^{\prime},\tau^{\prime},I^{\prime})$ be two para-Hermitian vector bundles.
  A bundle map $F:E \rightarrow E^{\prime}$ is called a para-holomorphic bundle map if $F \circ I = I^{\prime} \circ F$.
\end{dfn}


\begin{exam}\upshape\label{example of para-Kahler bundle}
  Consider the Whitney sum $E = TM \oplus T^{*}M$ on $M$, we define
  \begin{align*}
  \tau(X \oplus \alpha, Y \oplus \beta) &:= \beta(X) + \alpha(Y),\\
  I(X \oplus \alpha) &:= X \oplus (-\alpha)
\end{align*}
  for $X \oplus \alpha, Y \oplus \beta \in TM \oplus T^{*}M$.
  Then $(TM \oplus T^{*}M, \tau, I)$ is a para-Hermitian vector bundle.
\end{exam}

\vskip\baselineskip
\subsection{Quasi-Codazzi structure}
\label{sec:3-3}
Let $M$ be an $n$-dimensional $C^{\infty}$-manifold and $(E, \tau, I)$ a para-Hermitian vector bundle of $\rank 2n$ on $M$.
Set $\omega(X,Y):=\tau(X,IY)$ to be the symplectic form.

We consider a bundle map $\Phi:TM \rightarrow E$ such that, for every $p \in M$, $\rank\Phi_{p} = n$ and $\Imag \Phi_{p} = \Phi_{p}(T_{p}M)$ is a Lagrange subspace of $E_{p}$, i.e., $\omega$ vanishes on $\Imag \Phi_{p}$.
Then, $\Phi$ or $\Imag \Phi$ is a Lagrange subbundle.

\[
\xymatrix{
TM \ar[rd]_{\pi_{TM}} \ar[rr]^{\Phi} & & E \ar[ld]^{\pi_{E}} \\
& M &
}
\]

We call $E^{+}$ (resp. $E^{-}$) the {\em eigenbundle} of the eigenvalue $1$ (resp. $-1$).
For the above $\Phi$, we have two bundle maps $\Phi^{+}:=\pi^{+} \circ \Phi$ and $\Phi^{-}:=\pi^{-} \circ \Phi$, where $\pi^{+}:E \rightarrow E^{+}$ and $\pi^{-}:E \rightarrow E^{-}$ are projections with respect to $E=E^{+}\oplus E^{-}$.

Let $\nabla$ be a connection on $E^{+}$.

\begin{dfn}\upshape
  \label{dfn dual connection}
  The dual connection $\nabla^{*}$ on $E^{-}$ is defined by the equality
  \[
  X\tau(\nu^{+}, \nu^{-}) = \tau(\nabla_{X}\nu^{+}, \nu^{-}) + \tau(\nu^{+}, \nabla^{*}_{X}\nu^{-})
  \]
  for $X\in \Gamma(TM)$, $\nu^{+}\in \Gamma(E^{+})$ and $\nu^{-} \in \Gamma(E^{-})$.
\end{dfn}

In fact, $\nabla^{*}$ exists uniquely since $\tau$ is of $(n,n)$-type.
In particular, it holds that $(\nabla^{*})^{*} = \nabla$.
From now on, we use the notation
\[
\nabla^{+}:=\nabla \text{ on } E^{+},\quad \nabla^{-}:= \nabla^{*} \text{ on } E^{-}.
\]

\begin{dfn}\upshape
  We define a symmetric $(0,2)$-tensor $h$ on $M$ by the pullback of $\tau$.
  \[
  h(Y,Z) := \tau(\Phi(Y), \Phi(Z))
  \]
  with $Y,Z \in \Gamma(TM)$.
\end{dfn}

Below we use the following notation for $X,Y,Z\in \Gamma (TM)$:
\begin{gather*}
  \xi^{+}:=\Phi^{+}(X),\quad \xi^{-}:=\Phi^{-}(X),\quad \eta^{+} := \Phi^{+}(Y),\quad \eta^{-}:=\Phi^{-}(Y),\\ \zeta^{+}:=\Phi^{+}(Z), \quad \zeta^{-}:=\Phi^{-}(Z).
\end{gather*}

In particular, it holds that
\[
\tau(\eta^{+}, \zeta^{+}) = \tau(\eta^{-},\zeta^{-}) = 0
\]
and
\[
\tau(\eta^{+}\oplus\eta^{-},\zeta^{+}\oplus\zeta^{-}) = \tau(\eta^{+},\zeta^{-}) + \tau(\eta^{-},\zeta^{+}).
\]

\begin{lem}\upshape
  \label{lemm:tau is symmetric}
  It holds that
  \[
  \tau(\eta^{+},\zeta^{-}) = \tau(\eta^{-},\zeta^{+})=\frac{1}{2}h(Y,Z)
  \]
  for $Y,Z\in \Gamma(TM)$.
\end{lem}

\proof
The condition that $\Imag \Phi$ is Lagrangian is rewritten as
\begin{align*}
  0 = \omega(\Phi(Y), \Phi(Z)) &= \omega(\eta^{+} \oplus \eta^{-},\zeta^{+} \oplus \zeta^{-})  \\
  &= \tau(\eta^{+} \oplus \eta^{-}, I(\zeta^{+} \oplus \zeta^{-})) \\
  &= \tau(\eta^{+} \oplus \eta^{-}, \zeta^{+} \oplus -\zeta^{-}) \\
  &= -\tau(\eta^{+},\zeta^{-}) + \tau(\eta^{-}, \zeta^{+}).
\end{align*}
Also, we see that
\begin{align*}
  h(Y,Z) = \tau(\Phi(Y), \Phi(Z))
  &= \tau(\eta^{+},\zeta^{-}) + \tau(\eta^{-}, \zeta^{+}) \notag \\
  &= 2\tau(\eta^{+},\zeta^{-})\notag \\
  &=2\tau(\eta^{-},\zeta^{+}).
\end{align*}
\qed

\begin{lem}\upshape
  \label{h and Phi}
  The symmetric $(0,2)$-tensor $h$ is non-degenerate if and only if $\Phi^{+}$ and $\Phi^{-}$ are isomorphisms, i.e., they have $\rank n$ everywhere.
\end{lem}
\proof
Pick $p\in M$.
A matrix representation of Lemma \ref{lemm:tau is symmetric} is
\[
\frac{1}{2}{\mbox{\boldmath $h$}} = ({\mbox{\boldmath $\Phi$}}^{+})^{\top} {\mbox{\boldmath $\tau$}} {\mbox{\boldmath $\Phi$}}^{-}.
\]
Hence, since $\mbox{\boldmath $\tau$}$ is a regular matrix, $\mbox{\boldmath $h$}$ is a regular matrix if and only if $\mbox{\boldmath $\Phi$}^{+}$ and $\mbox{\boldmath $\Phi$}^{-}$ are regular matrices.
\qed

\begin{lem}\upshape
  \label{isom and dual}
  Let $h$ be non-degenerate, then two connections $\widetilde{\nabla}^{+}$ and $\widetilde{\nabla}^{-}$ of $TM$ are defined by
  \[
  \Phi^{+}(\widetilde{\nabla}^{+}_{X}Y):= \nabla^{+}_{X}\Phi^{+}(Y), \qquad
  \Phi^{-}(\widetilde{\nabla}^{-}_{X}Y):= \nabla^{-}_{X}\Phi^{-}(Y)
  \]
  for $X,Y\in \Gamma(TM)$, and they satisfy
  \[
  Xh(Y,Z) =  h(\widetilde{\nabla}^{+}_{X}Y, Z) + h(Y, \widetilde{\nabla}^{-}_{X}Z)
  \]
  for $X,Y,Z \in \Gamma(TM)$, i.e., $\til{\nabla}^{+}$ and $\til{\nabla}^{-}$ are mutually dual with respect to $h$.
\end{lem}

\proof
From Proposition \ref{dfn dual connection} and $h(Y,Z)=2\tau(\eta^{+},\zeta^{-})$,
\begin{align*}
  Xh(Y,Z) &= 2X\tau(\eta^{+},\zeta^{-}) \\
  &= 2\{ \tau(\nabla^{+}_{X}\eta^{+}, \zeta^{-}) + \tau(\eta^{+}, \nabla^{-}_{X}\zeta^{-})\} \\
  &= 2\tau(\Phi^{+}(\widetilde{\nabla}^{+}_{X}Y), \Phi^{-}(Z)) + 2\tau(\Phi^{+}(Y), \Phi^{-}(\widetilde{\nabla}^{-}_{X}Z)) \\
  &= h(\widetilde{\nabla}^{+}_{X}Y,  Z) + h(Y, \widetilde{\nabla}^{-}_{X}Z).
\end{align*}
\qed

\begin{dfn}\upshape
  \label{cubic-tensor}
  We define a $(0,3)$-tensor $C$ on $M$ by
  \[
  C(X,Y,Z) = -2\{ \tau(\nabla^{+}_{X}\eta^{+},\zeta^{-}) - \tau(\zeta^{+}, \nabla^{-}_{X}\eta^{-}) \}.
  \]
  Note that $C(X,Y,Z) = C(X,Z,Y)$, but $C$ is not totally symmetric in general.
\end{dfn}

\begin{lem}\upshape
  \label{cubic tensor if non-degenerate}
  If $h$ is non-degenerate, then $C$ is nothing but the Amari-Chentsov tensor $\til{\nabla}^{+}h$ in the original sense.
\end{lem}
\proof
We see that
\begin{align*}
C(X,Y,Z) &=  h(Z, \widetilde{\nabla}^{-}_{X}Y) - h(\widetilde{\nabla}^{+}_{X}Y,Z) \\
& = Xh(Y,Z) - h(\widetilde{\nabla}^{+}_{X}Y,Z) - h(\widetilde{\nabla}^{+}_{X}Z,Y) \\
&= (\widetilde{\nabla}^{+}h)(X,Y,Z)
\end{align*}
by using the definition of $\tau$ and the duality of $\widetilde{\nabla}^{+}$ and $\widetilde{\nabla}^{-}$.
\qed

\begin{dfn}\upshape
  \cite{HHNSUY, Kossowski}
  \label{Kossowski}
  The {\em Kossowski pseudo-connection} $\Gamma$ is a smooth map $\Gamma(TM) \times \Gamma(TM) \times \Gamma(TM) \rightarrow C^{\infty}(M)$ as
  \begin{align*}
    \Gamma(X,Y,Z) &= X\tau(\eta^{+}, \zeta^{-}) + Y\tau(\zeta^{+}, \xi^{-}) - Z\tau(\xi^{+},\eta^{-})\\
    &\qquad - \tau(\xi^{+}, \Phi^{-}([Y,Z])) + \tau(\eta^{+},\Phi^{-}([Z,X])) + \tau(\zeta^{+},\Phi^{-}([X,Y])).
  \end{align*}
\end{dfn}

\begin{lem}\upshape
  \label{Kossowski if non-degenerate}
  Assume that $h$ is non-degenerate.
  Let $\Gamma$ be the Kossowski pseudo-connection and $\nabla^{(0)}$ the Levi-Civita connection with respect to $h$.
  Then it holds that, for $X,Y,Z\in \Gamma(TM)$,
  \[
  h(\nabla^{(0)}_{X}Y,Z) = \Gamma(X,Y,Z).
  \]
\end{lem}

\proof
By definition of the Kossowski pseudo-connection,
\begin{align*}
  \Gamma(X,Y,Z) &= \{ X\tau(\eta^{+}, \zeta^{-}) + Y\tau(\zeta^{+}, \xi^{-}) - Z\tau(\xi^{+},\eta^{-})\\
  &\qquad - \tau(\xi^{+}, \Phi^{-}([Y,Z])) + \tau(\eta^{+},\Phi^{-}([Z,X])) + \tau(\zeta^{+},\Phi^{-}([X,Y])) \}\\
  &= \frac{1}{2} \{ Xh(Y,Z) + Yh(Z,X) - Zh(X,Y)\\
  &\qquad -h(X,[Y,Z])+h(Y,[X,Z])+h(Z,[X,Y]) \}\\
  &= h(\nabla^{(0)}_{X}Y,Z).
\end{align*}
The last equality is a well-known formula of Levi-Civita connection (see \cite[Proposition 2.3]{Kobayashi-Nomizu}).
\qed

\begin{dfn}\upshape
  \cite{Kobayashi-Nomizu}
  \label{curvature and relatively torsion}
  We define the {\em curvature} $R^{\nabla^{+}}$ with respect to $\nabla^{+}$ by
  \[
  R^{\nabla^{+}}(X,Y)\nu^{+} := (\nabla^{+}_{X} \nabla^{+}_{Y} - \nabla^{+}_{Y} \nabla^{+}_{X} - \nabla^{+}_{[X,Y]})\nu^{+}
  \]
  for $X,Y\in \Gamma(TM)$ and $\nu^{+} \in \Gamma(E^{+})$.
  In entirely the same way, we define $R^{\nabla^{-}}$ with respect to $\nabla^{-}$.
\end{dfn}

The metric $\tau$ and the curvature have the following relation.

\begin{prop}\upshape
  \label{curvature and metric}
  It holds that
  \[
  \tau(R^{\nabla^{+}}(X,Y)\nu^{+}, \nu^{-}) + \tau(\nu^{+}, R^{\nabla^{-}}(X,Y)\nu^{-}) = 0
  \]
  for $X,Y\in\Gamma(TM)$, $\nu^{+}\in\Gamma(E^{+})$, and $\nu^{-}\in \Gamma(E^{-})$.
\end{prop}
\proof
This is shown by a direct computation using Definition \ref{dfn dual connection} as follows:
\begin{align*}
  & \quad \tau(R^{\nabla^{+}}(X,Y)\nu^{+}, \nu^{-})\\
  &= \tau(\nabla^{+}_{X}\nabla^{+}_{Y}\nu^{+} - \nabla^{+}_{Y}\nabla^{+}_{X}\nu^{+} - \nabla^{+}_{[X,Y]}\nu^{+}, \nu^{-})\\
  &= \tau(\nabla^{+}_{X}\nabla^{+}_{Y}\nu^{+}, \nu^{-}) - \tau(\nabla^{+}_{Y}\nabla^{+}_{X}\nu^{+}, \nu^{-}) - \tau(\nabla^{+}_{[X,Y]}\nu^{+}, \nu^{-})\\
  &= \{ X\tau(\nabla^{+}_{Y}\nu^{+},\nu^{-}) - \tau(\nabla^{+}_{Y}\nu^{+},\nabla^{-}_{X}\nu^{-}) \} -\{ Y\tau(\nabla^{+}_{X}\nu^{+},\nu^{-}) - \tau(\nabla^{+}_{X}\nu^{+}, \nabla^{-}_{Y}\nu^{-}) \}\\
  &\qquad \qquad -\{ [X,Y] \tau(\nu^{+},\nu^{-}) - \tau(\nu^{+}, \nabla^{-}_{[X,Y]} \nu^{-}) \} \\
  &= X \left\{ Y\tau(\nu^{+}, \nu^{-}) - \tau(\nu^{+}, \nabla^{-}_{Y}\nu^{-}) \right\} - \{Y\tau(\nu^{+},\nabla^{-}_{X}\nu^{-}) - \tau(\nu^{+},\nabla^{-}_{Y} \nabla^{-}_{X}\nu^{-}) \}\\
  &\qquad \qquad -Y\{ X\tau(\nu^{+},\nu^{-}) - \tau(\nu^{+}, \nabla^{-}_{X}\nu^{-}) \} + \{ X\tau(\nu^{+},\nabla^{-}_{Y}\nu^{-}) - \tau(\nu^{+}, \nabla^{-}_{X}\nabla^{-}_{Y}\nu^{-}) \}\\
  &\qquad \qquad \qquad -[X,Y]\tau(\nu^{+},\nu^{-}) + \tau(\nu^{+}, \nabla^{-}_{[X,Y]}\nu^{-})\\
  &= -\tau(\nu^{+}, \nabla^{-}_{X}\nabla^{-}_{Y}\nu^{-} - \nabla^{-}_{Y}\nabla^{-}_{X}\nu^{-} - \nabla^{-}_{[X,Y]}\nu^{-})\\
  &= - \tau(\nu^{+}, R^{\nabla^{-}}(X,Y)\nu^{-}).
\end{align*}
\qed

\begin{cor}\upshape
  If $R^{\nabla^{+}}\equiv 0$, then $R^{\nabla^{-}} \equiv 0$, and vice-versa.
  That is, $\nabla^{+}$ is flat if and only if $\nabla^{-}$ is flat (we say that $\nabla^{+}$ and $\nabla^{-}$ are dually flat).
\end{cor}
\proof
It immediately follows from Proposition \ref{curvature and metric}.
\qed

\vskip\baselineskip
Now we seek for a singular version of Proposition \ref{prop of Codazzi str}.
Consider the following conditions:
\begin{itemize}
  \item[(i)]
  $\nabla^{+}$ is relatively torsion-free;
  \item[(ii)]
  $\nabla^{-}$ is relatively torsion-free;
  \item[(iii)]
  $C$ is totally symmetric, i.e., it is verified by
  \[
  C(X,Y,Z) = C(Y,X,Z);
  \]
  \item[(iv)]
  It holds that
  \[
  \Gamma(X,Y,Z)=\tau(\nabla^{+}_{X}\eta^+,\zeta^{-}) + \tau(\zeta^{+},\nabla^{-}_{X}\eta^{-}).
  \]
\end{itemize}

Notice that if $h$ is non-degenerate, then the above four conditions (i)-(iv) coincide  with the conditions in Proposition \ref{prop of Codazzi str}, respectively, that follows from Lemma \ref{isom and dual}, Lemma \ref{cubic tensor if non-degenerate} and Lemma \ref{Kossowski if non-degenerate}.

We present a generalization of Proposition \ref{prop of Codazzi str} to a singular setting as follows.

\begin{prop}\upshape
  \label{prop of quasi-Codazzi str}
  $ $
  \begin{itemize}
    \item[(1)]
    If (i) and (ii) are satisfied, then (iii) and (iv) is satisfied;
    \item[(2)]
    If $\Phi^{+}$ is an isomorphism and (i) and either (iii) or (iv) are satisfied, then (ii) is satisfied;
    \item[(3)]
    If $\Phi^{-}$ is an isomorphism and (ii) and either (iii) or (iv) are satisfied, then (i) is satisfied;
    \item[(4)]
    If $h$ is non-degenerate and two of the there conditions are satisfied, then remaining are satisfied.
  \end{itemize}
\end{prop}

\proof
By Definition \ref{def of relative torsion and coherent tangent bundle}, Lemma \ref{lemm:tau is symmetric} and Definition \ref{cubic-tensor},
\begin{align*}
  -\frac{1}{2} &\{ C(X,Y,Z) - C(Y,X,Z) \}\\
  &= \{ \tau(\nabla^{+}_{X}\eta^{+},\zeta^{-}) - \tau(\zeta^{+},\nabla^{-}_{X}\eta^{-}) \} - \{\tau(\nabla^{+}_{Y}\xi^{+},\zeta^{-}) - \tau(\zeta^{+},\nabla^{-}_{Y}\xi^{-}) \}\\
  &=\tau(T^{\nabla^{+}}(X,Y), \zeta^{-}) + \tau(\Phi^{+}([X,Y]),\zeta^{-})\\
  &\qquad \qquad \qquad- \tau(\zeta^{+}, T^{\nabla^{-}}(X,Y)) -\tau(\zeta^{+}, \Phi^{-}([X,Y]))\\
  &=\tau(T^{\nabla^{+}}(X,Y),\zeta^{-}) - \tau(\zeta^{+},T^{\nabla^{-}}(X,Y)).
\end{align*}
Hence we get
\[
(\text{iii}) \Longleftrightarrow \tau(T^{\nabla^{+}}(X,Y),\zeta^{-}) = \tau(\zeta^{+},T^{\nabla^{-}}(X,Y)).
\]
Also, we can easily check
\[
(\text{iv}) \Longleftrightarrow \tau(T^{\nabla^{+}}(X,Y),\zeta^{-}) = -\tau(\zeta^{+},T^{\nabla^{-}}(X,Y)).
\]
Hence, (1)-(4) are immediately shown.
This complete the proof.
\qed

\begin{dfn}\upshape
  \label{def of quasi-Codazzi structure}
  Let $(h,(E,\tau,I),\Phi,\nabla^{+},\nabla^{-})$ be as above, i.e., $(E,\tau,I)$ is a para-Hermitian vector bundle with a Lagrange subbundle $\Imag \Phi$ and connections $\nabla^{\pm}$ and a symmetric $(0,2)$-tensor $h=\Phi^{*}\tau$.
  We call it a {\em quasi-Codazzi structure} if $\nabla^{+}$ and $\nabla^{-}$ are relatively torsion-free, i.e., $\Phi^{+}$ and $\Phi^{-}$ define two coherent tangent bundles in the sense of Definition \ref{def of relative torsion and coherent tangent bundle}, and a manifold $M$ with a quasi-Codazzi structure is called a {\em quasi-Codazzi manifold}.
  Then we call the symmetric cubic tensor $C$ the {\em generalized Amari-Chentsov tensor}.
\end{dfn}
\[
\xymatrix{
E^{+} \ar[rd]_{\pi_{E^{+}}} & TM \ar[l]_{\Phi^{+}} \ar[r]^{\Phi^{-}} \ar[d]_{\pi_{TM}} & E^{-} \ar[ld]^{\pi_{E^{-}}} \\
& M &
}
\]

A Codazzi structure in ordinary sense defines a quasi-Codazzi structure in a canonical way as follows.

\begin{exam}\upshape
  \label{canonical quasi-Codazzi}
  Suppose that we are given a Codazzi structure $(h,\nabla,\nabla^{*})$ on $M$.
  Let $(TM\oplus T^{*}M, \tau, I)$ be a para-Hermitian vector bundle as Example \ref{example of para-Kahler bundle}.
  Set $\Xi:TM\rightarrow TM \oplus T^{*}M$ by
  \[
  \Xi(X):= \left( id(X), \frac{1}{2}\til{h}(X) \right)= \left( X,\frac{1}{2} h(X,-) \right)
  \]
  and connections $\nabla^{\pm}$ by
  \[
  \nabla^{+}_{X}Y= \nabla_{X}Y,\qquad \nabla^{-}_{X} \til{h}(Y)= \til{h}( \nabla^{*}_{X}Y).
  \]
  Then $\Xi$ is Lagrangian and $\nabla^{\pm}$ are relatively torsion-free.
  Thus,
  \[
  (h,(TM\oplus T^{*}M,\tau,I),\Xi, \nabla^{+},\nabla^{-})
  \]
  is a quasi-Codazzi structure.
\end{exam}

In Proposition \ref{prop:quasi-Codazzi structure with non-degenerate => Codazzi structure} below, we show that a quasi-Codazzi structure with non-degenerate $h$ is equivalent to the Codazzi structure given in Example \ref{canonical quasi-Codazzi}.

\vskip\baselineskip
\subsection{Example}\label{sec:3-3.5}
We present examples of quasi-Codazzi structure.
\vskip\baselineskip
\noindent
(1)First, let $(N,\tau,I)$ be a para-Hermitian manifold with $TN = E^{+}\oplus E^{-}$, the decomposition to eigenspaces of $I$, and let $\nabla = \nabla^{+}\oplus \nabla^{-}$ be an invariant connection with respect to $
I$.
Then any Lagrange submanifold $M$ of $N$ admits a quasi-Codazzi structure with $\Phi:TM\rightarrow TN$ and $h:=\Phi^{*}\tau$ which can be degenerate in general.
That may possess a {\em Frobenius-like structure} \cite[\S5.1]{Nakajima-Ohmoto}.
Lagrange submanifolds in para-K\"ahler manifolds would be also meaningful in the context of quasi-Codazzi structure.

\vskip\baselineskip
\noindent
(2)The second example is {\em front bundle} introduced by Saji-Umehara-Yamada \cite{SUY-F2, SUY-F1}.
A full discussion will be given in the forthcoming paper of the author \cite{Kayo}.
A typical case arises from a singular wavefront in an ambient Riemannian manifold with constant curvature \cite{SUY-F2}.
That is explained shortly below.

  Let $M$ be an $n$-dimensional manifold, $(N,g)$ a complete Riemannian $(n+1)$-dimensional manifold of constant curvature $c$, and $\nabla$ the Levi-Civita connection on $N$ with respect to $g$.
  A $C^{\infty}$-map $f:M\rightarrow N$ is called a {\em front} or {\em wavefront} if for each $p\in M$, there exists a neighborhood $U$ of $p$ and a unit vector field $\nu$ which satisfies $g(df(X),\nu)=0$ for $X\in\Gamma(TU)$, and the map $\nu:U\ni q \mapsto \nu_{q}\in T_{1}N$ is an immersion as a smooth section of $T_{1}N$, where $T_{1}N$ is the unit tangent bundle of $N$.
  Assume that $f:M\rightarrow N$ is a wavefront, we get a coherent tangent bundle.
  In fact, let $\mathcal{E}$ be defined as $\mathcal{E}:=\{ \xi \in f^{*}TN | \langle \xi,\nu \rangle=0 \}$, and $\langle \;,\;\rangle$ a metric on $\mathcal{E}$ induced from $g$, and $D$ a connection on $\mathcal{E}$ defined by taking the tangential part of $\nabla$, and $\phi:=df$.
  Then $(\mathcal{E}, \phi,D)$ be a coherent tangent bundle.
  Moreover, by the condition of constant curvature, we can get another coherent tangent bundle $(\mathcal{E},\psi,D)$ with
  \[
  \psi_{p}:T_{p}M \ni X_{p} \mapsto \nabla_{X_{p}}\nu_{p} \in \mathcal{E}_{p}
  \]
  for each $p\in M$.
  This defines a quasi-Codazzi structure on $M$ with $E^{+} =E^{-}=\mathcal{E}$, $\Phi^{+} = \phi$, $\Phi^{-}=\psi$ : $\Phi=(\phi,\psi)$ defines a Lagrange subbundle of $\mathcal{E}\oplus \mathcal{E}$, that is implicitly described in \cite[Example 2.2]{SUY-F2}.

\vskip\baselineskip
\subsection{Isomorphism of quasi-Codazzi structures}\label{sec:3-4}

In this section, we introduce an equivalence relation among quasi-Codazzi structures on $M$.

%
%
%
For $i=1,2$, let
\[
A_{i} := (h_{i},(E_{i},\tau_{i},I_{i}),\Phi_{i},\nabla^{+}_{i},\nabla^{-}_{i})
\]
be quasi-Codazzi structures on $M$.

\begin{dfn}\upshape\label{def of isomorphic}
  Quasi-Codazzi structures $A_{1}$ and $A_{2}$ are isomorphic if there exists a bundle isomorphism $F:E_{1} \rightarrow E_{2}$ on $M$ which satisfies
  \begin{itemize}
    \item [(i)]
    $F \circ \Phi_{1} = \Phi_{2}$;
    \item[(ii)]
    $F \circ I_{1} = I_{2}  \circ F$;
    \item[(iii)]
    $\tau_{1}(\nu,\kappa) = \tau_{2}(F(\nu),F(\kappa))$ for $\nu,\kappa\in \Gamma(E_{1})$;
    \item[(iv)]
    $F(\nabla^{+}_{1\;X}\nu^{+}) = \nabla^{+}_{2\;X}F(\nu^{+})$ for $\nu^{+}\in \Gamma(E_{1}^{+})$ and $F(\nabla^{-}_{1\;X}\nu^{-}) = \nabla^{-}_{2\;X}F(\nu^{-})$ for $\nu^{-}\in \Gamma(E_{1}^{-})$.
  \end{itemize}
\end{dfn}

\begin{rem}\upshape
  Note that (i) and (ii) are satisfied if and only if the following diagram commutes:
  \begin{align*}
    \vcenter{
      \xymatrix{
        E^{+}_{1} \ar[d]_{F|_{E^{+}_{1}}} & TM \ar[l]_{\Phi^{+}_{1}} \ar[r]^{\Phi^{-}_{1}} \ar[d]_{id} & E^{-}_{1} \ar[d]^{F|_{E_{1}^{-}}} \\
        E_{2}^{+} & TM \ar[l]^{\Phi^{+}_{2}} \ar[r]_{\Phi^{-}_{2}} & E^{-}_{2}
        }
    }
  \end{align*}
\end{rem}

\begin{prop}\upshape\label{prop of isomorphic}
  Assume that $A_{1}$ and $A_{2}$ are isomorphic by $F:E_{1} \rightarrow E_{2}$.
  Then we have, for $X,Y,Z \in \Gamma(TM)$, $\nu^{+} \in \Gamma(E^{+})$, and $\nu^{-} \in \Gamma(E^{-})$, the following hold:
  \begin{itemize}
    \item[$(1)$]
    $F(R^{\nabla^{+}_{1}}(X,Y)\nu^{+}) = R^{\nabla^{+}_{2}}(X,Y)F(\nu^{+})$;
    \item[$(2)$]
    $F(R^{\nabla^{-}_{1}}(X,Y)\nu^{-}) = R^{\nabla^{-}_{2}}(X,Y)F(\nu^{-})$;
    \item[$(3)$]
    $h_{1}(X,Y) = h_{2}(X,Y)$.
    \item[$(4)$]
    $C_{1}(X,Y,Z) = C_{2}(X,Y,Z)$, where $C_{i}$ is the cubic tensor defined by Definition \ref{cubic-tensor} with respect to $A_{i}$ for $i=1,2$.
  \end{itemize}
\end{prop}

\proof
This is shown by a straightforward calculation.
\qed

\begin{prop}\upshape
  \label{prop:quasi-Codazzi structure with non-degenerate => Codazzi structure}
  Let $(h,(E,\tau,I),\Phi,\nabla^{+},\nabla^{-})$ be a quasi-Codazzi structure.
  Suppose $\Phi^{+}$ and $\Phi^{-}$ are isomorphisms (i.e., $h$ is non-degenerate), and set connections $\til{\nabla}^{\pm}$ on $TM$ as in Lemma \ref{isom and dual}.
  Then this structure is isomorphic to the quasi-Codazzi structure $(h,(TM \oplus T^{*}M, \widehat{\tau}, \widehat{I}), id_{TM} \oplus (1/2)\til{h}, \widehat{\nabla}^{+}, \widehat{\nabla}^{-})$ in Example \ref{canonical quasi-Codazzi}.
\end{prop}

\proof
Set $F:=(\Phi^{+})^{-1} \oplus (1/2)\til{h} \circ (\Phi^{-})^{-1}$.
Then $F \circ \Phi = id_{TM} \oplus (1/2)\til{h}$, $F \circ I = \widehat{I} \circ F$, and
\begin{align*}
\tau(\Phi(Y),\Phi(Z)) &= h(Y,Z)\\
&= \frac{1}{2}h(Y,Z) + \frac{1}{2}h(Y,Z)\\
&= \widehat{\tau} \left( Y,\frac{1}{2}\til{h}(Z) \right) + \widehat{\tau} \left( \frac{1}{2}\til{h}(Y), Z \right)\\
&= \widehat{\tau} \left( Y \oplus \frac{1}{2}\til{h}(Y), Z \oplus \frac{1}{2}\til{h}(Z) \right)\\
&= \widehat{\tau}(F \circ \Phi(Y), F \circ \Phi(Z)).
\end{align*}
Furthermore, for $X,Y\in \Gamma(TM)$,
\begin{align*}
  F(\nabla^{+}_{X}\eta^{+}) = (\Phi^{+})^{-1}  \circ \Phi^{+}(\til{\nabla}^{+}_{X}Y)=\til{\nabla}^{+}_{X}Y = \widehat{\nabla}^{+}_{X}Y = \widehat{\nabla}^{+}_{X}F(\eta^{+})
\end{align*}
and
\begin{align*}
  F(\nabla^{-}_{X}\eta^{-}) = \frac{1}{2}\til{h} \circ (\Phi^{-})^{-1}  \circ \Phi^{-}(\til{\nabla}^{-}_{X}Y)=\frac{1}{2}\til{h}(\til{\nabla}^{-}_{X}Y) = \widehat{\nabla}^{-}_{X} \frac{1}{2}\til{h}(Y) =  \widehat{\nabla}^{-}_{X} F(\eta^{-}).
\end{align*}
\qed

\vskip\baselineskip
\subsection{Weak contrast functions}
\label{sec:3-5}
Given a Codazzi structure $(h,\nabla,\nabla^{*})$, we find a contrast function which induces $h$, $C$, $\nabla$, and $\nabla^{*}$ \cite{Eguchi}.
Let $(h,(E,\tau,I),\Phi,\nabla^{+},\nabla^{-})$ be a quasi-Codazzi structure on $M$.
Since this structure has a possibly degenerate metric $h$ and a totally symmetric cubic-tensor $C$, a weak contrast function $\rho$ exists such that
\begin{align*}
  h(X,Y) &= -\rho[X|Y], \\
  C(X,Y,Z) &= -\rho[Z|XY] + \rho[XY|Z]
\end{align*}
by Theorem \ref{weak contrast function}.
This weak contrast function $\rho$ and the quasi-Codazzi structure have the following relation.

\begin{prop}\upshape\label{relation of a weak contrast 1}
  (cf. \cite[Lemma 2.1]{Matumoto})
  There exist relations, for $X,Y,Z \in \Gamma(TM)$,
  \begin{itemize}
    \item [(1)]
    $-\rho[XY|Z]=  \Gamma(X,Y,Z) - \dfrac{1}{2} C(X,Y,Z)$;
    \item [(2)]
    $-\rho[Z|XY] = \Gamma(X,Y,Z) + \dfrac{1}{2}C(X,Y,Z)$.
  \end{itemize}
\end{prop}
\proof
The following equations are satisfied by properties of $\rho$,
\[
\begin{cases}
  X\rho[Y|Z] = \rho[XY|Z] + \rho[Y|XZ];\\
  \rho[X|[Y,Z]] = \rho[X|YZ] -\rho[X|ZY].
\end{cases}
\]
Then,
\begin{align*}
  2\Gamma(X,Y,Z) &= Xh(Y,Z) + Yh(Z,X) -Zh(X,Y)\\
  & \qquad \qquad -h(X,[Y,Z]) + h(Y,[Z,X]) + h(Z,[X,Y])\\
  &= -X\rho[Y|Z] -Y\rho[Z|X]+Z\rho[X|Y]\\
  &   \qquad \qquad +\rho[X|[Y,Z]]-\rho[Y|[Z,X]]-\rho[Z|[X,Y]]\\
  &= -\{ \rho[XY|Z] + \rho[Y|XZ] \}  - \{ \rho[YZ|X] + \rho[Z|YX] \}\\
  &  \quad  +\{ \rho[ZX|Y] + \rho[X|ZY] \} + \{ \rho[X|YZ] - \rho[X|ZY] \}\\
  &  \qquad - \{ \rho[Y|ZX] -  \rho[Y|XZ] \} -  \{ \rho[Z|XY] - \rho[Z|YX] \}\\
  &= -\rho[XY|Z] - \rho[YZ|X]+\rho[ZX|Y]\\
  & \qquad \qquad +\rho[X|YZ] -\rho[Y|ZX]- \rho[Z|XY]\\
  &=-\rho[XY|Z] -C(Y,Z,X) + C(Z,X,Y) - \rho[Z|XY]\\
  &=-\rho[XY|Z] - \rho[Z|XY].
\end{align*}
The last equality comes from that $C$ is totally symmetric.
Hence we get equations
\[
2\Gamma(X,Y,Z) - C(X,Y,Z) = -2\rho[XY|Z]
\]
and
\[
2\Gamma(X,Y,Z) + C(X,Y,Z) = -2\rho[Z|XY].
\]
\qed

\begin{cor}\upshape\label{relation of a weak contrast}
  The following equations are satisfied.
  \begin{itemize}
    \item [(1)]
    $\tau(\nabla^{+}_{X}\eta^{+}, \zeta^{-}) = -\dfrac{1}{2}\rho[XY|Z]$;
    \item[(2)]
    $\tau(\zeta^{+}, \nabla^{-}_{X}\eta^{-}) = -\dfrac{1}{2}\rho[Z|XY]$.
  \end{itemize}
\end{cor}

\proof
The quasi-Codazzi structure implies that
\[
\Gamma(X,Y,Z) = \tau(\nabla^{+}_{X}\eta^{+},\zeta^{-}) + \tau(\zeta^{+},\nabla^{-}_{X}\eta^{-}).
\]
Therefore, by Proposition \ref{relation of a weak contrast 1},
\begin{align*}
  -\rho[XY|Z] &= \Gamma(X,Y,Z)-\frac{1}{2}C(X,Y,Z)\\
  &= \{\tau(\nabla^{+}_{X}\eta^{+},\zeta^{-}) + \tau(\zeta^{+},\nabla^{-}_{X}\eta^{-})\} + \{ \tau(\nabla^{+}_{X}\eta^{+},\zeta^{-}) - \tau(\zeta^{+}, \nabla^{-}_{X} \eta^{-}) \}\\
  &= 2\tau(\nabla^{+}_{X}\eta^{+}, \zeta^{-}).
\end{align*}
We get (2) in the same way.
\qed

\begin{cor}\upshape
  The following hold, for $W,X,Y,Z \in \Gamma(TM)$,
  \begin{itemize}
      \item[$(1)$]
      $\displaystyle \tau(R^{\nabla^{+}}(X,Y)\zeta^{+},\Phi^{-}(W)) = \tau(\nabla^{+}_{X}\nabla^{+}_{Y}\zeta^{+},\Phi^{-}(W)) - \tau(\nabla^{+}_{Y}\nabla^{+}_{X}\zeta^{+},\Phi^{-}(W))$ \\
      $\qquad \qquad \qquad \qquad \qquad \qquad \qquad \qquad + \dfrac{1}{2}\rho[XYZ|W] - \dfrac{1}{2}\rho[YXZ|W]$;
      \item[$(2)$]
      $\displaystyle \tau(\Phi^{+}(W),R^{\nabla^{-}}(X,Y)\zeta^{-}) = \tau(\Phi^{+}(W), \nabla^{-}_{X}\nabla^{-}_{Y}\zeta^{-}) - \tau(\Phi^{+}(W),\nabla^{-}_{Y}\nabla^{-}_{X}\zeta^{-})$ \\
      $\qquad \qquad \qquad \qquad \qquad \qquad \qquad \qquad + \dfrac{1}{2}\rho[W|XYZ] - \dfrac{1}{2}\rho[W|YXZ].$
  \end{itemize}
\end{cor}

\proof
We show these equations by definition of the curvature and Corollary \ref{relation of a weak contrast}.

\qed

\vskip\baselineskip
\vskip\baselineskip
\section{Relations between quasi-Hessian structure and quasi-Codazzi structure}
\label{sec:4}

\vskip\baselineskip
\subsection{Flat connections and Poincar\'{e} Lemma}
\label{sec:4-2}
Let $M$ be an $n$-dimensional manifold, $\E \rightarrow M$ a vector bundle of $\rank r$ ($n \leq r$) on $M$, $\nabla$ a connection on $\E$, and $\Phi:TM \rightarrow \E$ a bundle map which is injective on each fiber, i.e., $\Phi_{p}:T_{p}M \rightarrow \E_{p}$ has $\rank n$.

%
%

Suppose that a connection $\nabla$ is curvature-free and relatively torsion-free with respect to $\Phi$.
As well known (e.g., \cite{Kobayashi-Nomizu}), for $\nabla$ is curvature-free, we can take a coordinate neighborhood $(U,\phi = (u_{1}, \cdots, u_{n}))$ and a local frame $(e_{1}, \cdots, e_{r})$ of $\E$ on $U$ such that $\nabla e_{i}=0$ $(1 \leq i \leq r)$.
We may assume that $\phi(U)$ is an open ball in $\R^{n}$.

We write $\Phi(\frac{\rd}{\rd u_{i}})=\sum_{j=1}^{r} \Phi_{ji}e_{j}$ on $U$, then
\begin{align*}
  \nabla_{\frac{\rd}{\rd u_{i}}} \Phi \left( \frac{\rd}{\rd u_{j}} \right) = \nabla_{\frac{\rd}{\rd u_{i}}} \sum_{k=1}^{r} \Phi_{kj} e_{k} &= \sum_{k=1}^{r} \left( \frac{\rd}{\rd u_{i}} \Phi_{kj} \right) e_{k} + \sum_{k=1}^{r}\Phi_{kj} \nabla_{\frac{\rd}{\rd u_{i}}}e_{k}\\
  &= \sum_{k=1}^{r} \left( \frac{\rd}{\rd u_{i}} \Phi_{kj} \right) e_{k}.
\end{align*}
The last equality holds by $\nabla e_{k}=0$.
Also,
\[
\nabla_{\frac{\rd}{\rd u_{j}}} \Phi \left( \frac{\rd}{\rd u_{i}} \right) = \sum_{k=1}^{r} \left( \frac{\rd}{\rd u_{j}} \Phi_{ki} \right)e_{k}.
\]
Since $\nabla$ is relatively torsion-free and $[ \frac{\rd}{\rd u_{i}}, \frac{\rd}{\rd u_{j}} ]=0$, we have
\begin{align}
  0 &= \nabla_{\frac{\rd}{\rd u_{i}}} \Phi \left( \frac{\rd}{\rd u_{j}} \right) -\nabla_{\frac{\rd}{\rd u_{j}}}\Phi \left( \frac{\rd}{\rd u_{i}} \right)\notag \\ 
  &= \sum_{k=1}^{r} \left( \frac{\rd}{\rd u_{i}} \Phi_{kj} - \frac{\rd}{\rd u_{j}} \Phi_{ki} \right)e_{k}.
  \label{integrability condition}
\end{align}

Consider the system of differential equations
\[
\frac{\rd}{\rd u_{1}}f_{k} = \Phi_{k1}, \cdots, \frac{\rd}{\rd u_{n}}f_{k} = \Phi_{kn}, \quad (f_
{k}:\phi(U)\rightarrow \R, (k=1,\cdots  r)),
\]
where $f_{k}$ $(1 \leq k \leq r)$ are $C^{\infty}$-functions on the open ball $\phi(U)\subset \R^{n}$.
By the Poincar\'e Lemma (e.g., \cite{Bott-Tu}), the integrability conditions of them are actually
\[
\frac{\rd}{\rd u_{i}} \Phi_{kj} = \frac{\rd}{\rd u_{j}} \Phi_{ki}
\]
for every $i,j \in \{ 1,\cdots,n \}$, that is \eqref{integrability condition}.
Therefore there exists $f=(f_{1}, \cdots, f_{r})^{\top} : \phi(U) \rightarrow \R^{r}$ such that, for $p\in U$,
\[
d(f \circ \phi)_{p} = \Phi_{p}:T_{p}U \rightarrow \E_{p}=\R^{r}.
\]
In particular, since $\Phi_{p}$ is injective for every $p\in U$, $f\circ \phi$ is a local embedding.

We summarize the above discussion.

\begin{prop}\upshape
  \label{prop of locally embedding}
  If $\nabla$ is curvature-free and relatively torsion-free with respect to $\Phi$, for any coordinate neighborhood $U$ of $M$ such that $U$ is contractible and $\E|_{U} \simeq U \times \R^{r}$,
  then there exists $f:U\rightarrow \R^{r}$ such that
  \[
  df\left( \frac{\rd}{\rd u_{i}}  \right) = pr \circ \Phi \left( \frac{\rd}{\rd u_{i}} \right) \qquad (1\leq i \leq n),
  \]
  where $pr:\E|_{U} \rightarrow \R^{r}$ is the projection.
\end{prop}

\vskip\baselineskip
\subsection{Quasi-Hessian structure}
\label{sec:4-1}
Consider the standard contact manifold $\R^{2n+1}=T^{*}\R^{n} \times \R$ with the Darboux coordinates $(\bx, \bp, z)$ \cite{AGV}.
It has a double Legendre (linear) fibration:
\[
\xymatrix{
\R^{n} _{\bbx} \times \R & \R^{2n+1} \ar[l] \ar[r] & \R^{n}_{\bbp} \times \R.
}
\]

\begin{dfn}\upshape
  \label{def of affine Legendre}
  An {\em affine Legendre equivalence} is an affine transformation $\Lcal :\R^{2n+1} \rightarrow \R^{2n+1}$ of the form
  \[
  \Lcal (\bx, \bp, z) = (\Lcal_{1}(\bx,\bp),\varphi(\bx,z)) = (A \bx+\bb, A^{\prime}\bp +\bb^{\prime}, z+\bc^{\top} \bx + d),
  \]
  with $\bb, \bb^{\prime} \in \R^{n},\; d\in \R, \;A,A^{\prime} \in GL_n(\R)$ satisfying that $A^{\prime}=(A^{\top})^{-1}$, and $\bb^{\prime} = A^{\prime}\bc$.
  Also, we call $\Lcal_{1}:\R^{2n}\rightarrow\R^{2n}$ an {\em affine Lagrange equivalence}.
\end{dfn}

\begin{dfn}\upshape
  \label{def of quasi-Hessian}
  A quasi-Hessian manifold is an affine manifold made up by gluing several Legendre submanifolds $L_{\alpha} \subset \R^{2n+1}$ via affine Legendre equivalences.
\end{dfn}

For more detail, see \cite[\S 3]{Nakajima-Ohmoto}.
Obviously, any quasi-Hessian manifold canonically admits a quasi-Codazzi structure.
In fact, locally, the differential of the projections of a Legendre submanifold $L_{\alpha}$ to the base and the fiber gives
\[
\xymatrix{
T\R^{n}_{\bbx} & TL_{\alpha} \ar[l]_{\Phi^{+}} \ar[r]^-{\Phi^{-}}  & T\R^{n}_{\bbp} = T^{*}\R^{n}_{\bbx}
}
\]
with a flat connection $\nabla$ on $\R^{n}_{\bbx}$ and its dual $\nabla^{*}$ on $\R^{n}_{\bbp}$, and these are compatible with affine Legendre equivalences.
Thus we are interested in the inverse, i.e., whether or not a quasi-Codazzi structure with flat connections defines a quasi-Hessian manifold.

\begin{thm}\upshape
  \label{thm of quasi-Hessian and quasi-Codazzi}
  Let $(h,(E=E^{+}\oplus E^{-},\tau,I),\Phi=\Phi^{+}\oplus\Phi^{-},\nabla^{+},\nabla^{-})$ be a quasi-Codazzi structure on $M$.
  If $\nabla^{+}$ and $\nabla^{-}$ are flat, then this structure naturally yields a quasi-Hessian structure on $M$.
\end{thm}

\proof
Let $U$ be a coordinate neighborhood which trivializes $E$.
Since $\nabla^{\pm}$ is a curvature-free connection on $E^{\pm}$, we get a frame $\{ e^{\pm}_{i} \}$ on $U$ which satisfy $\nabla^{\pm} e^{\pm}_{i}=0$ for every $i\in \{ 1,\cdots,n \}$.
By Definition \ref{dfn dual connection}, for $X\in \Gamma(TM), i,j\in \{ 1,\cdots,n \}$,
\[
X\tau(e_{i}^{+},e_{j}^{-}) = \tau(\nabla^{+}_{X}e_{i}^{+},e_{j}^{-}) + \tau(e_{i}^{+},\nabla^{-}_{X}e_{j}^{-}) = 0.
\]
Therefore $\tau(e_{i}^{+},e_{j}^{-})$ is a constant.
A matrix representation of $\tau$ is
\[
\tau=
\begin{pmatrix}
  \bzero & (a_{ij})_{n\times n}\\
  (a_{ij})^{\top}_{n\times n} & \bzero
\end{pmatrix},
\]
with $a_{ij} := \tau(e^{+}_{i},e^{-}_{j})$.
Since $\tau$ is non-degenerate, $(a_{ij})_{n\times n}$ is non-degenerate.
Set a frame $\{ \til{e}_{i}^{+} := \sum_{j=1}^{n} a^{ij}e_{j}^{+} \}$, where a matrix $(a^{ij})_{n\times n}$ is an inverse matrix of $(a_{ij})_{n\times n}$.
This frame satisfies $\nabla^{+} \til{e}_{i}^{+} = 0$ for $\nabla^{+}e_{j}^{+}=0$, and $\tau(\til{e}_{i}^{+},e_{j}^{-})=\delta_{ij}$, where $\delta_{ij}$ is the Kronecker symbol.
Hence we can retake $\{ e_{i}^{+} \}$ and $\{ e_{i}^{-} \}$ which satisfy, for every $i,j\in\{ 1,\cdots,n \}$,
\[
\nabla^{+}e_{i}^{+}=0, \quad \nabla^{-}e_{i}^{-}=0  \quad \text{and} \quad \tau(e_{i}^{+},e_{j}^{-})=\delta_{ij}.
\]
Set a connection $\nabla^{+}\oplus \nabla^{-}$ on $E$ which is defined by, for $X\in \Gamma(TM)$ and $\nu^{+}\oplus \nu^{-}\in \Gamma(E)$,
\[
(\nabla^{+}\oplus \nabla^{-})_{X}(\nu^{+}\oplus \nu^{-}):= (\nabla^{+}_{X}\nu^{+}) \oplus (\nabla^{-}_{X}\nu^{-}).
\]
Clearly, $\nabla^{+}\oplus \nabla^{-}$ is curvature-free and relatively torsion-free with respect to $\Phi$.
Thus, by Proposition \ref{prop of locally embedding}, we get an embedding
\[
f = f^{+}\oplus f^{-}:U \rightarrow (\R^{n}_{\bbx}\times\R^{n}_{\bbp}, (x_{1},\cdots,x_{n},p_{1},\cdots,p_{n}))
\]
with $\bx = (x_{1},\cdots,x_{n})$, $\bp=(p_{1},\cdots,p_{n})$ which satisfies
\[
df = df^{+}\oplus df^{-} = \Phi^{+}\oplus\Phi^{-}
\]
and
\[
\frac{\rd}{\rd x_{i}} = e_{i}^{+}, \quad \frac{\rd}{\rd p_{i}} = e_{i}^{-}
\]
for $i\in \{ 1,\cdots,n \}$.
Since the image of $\Phi$ is Lagrangian by the definition, $U$ can be realized as a Lagrange submanifold of $\R^{2n}=\R^{n}_{\bbx} \times \R^{n}_{\bbp}$.

Now, take another neighborhood $V$ with $U\cap V\neq \emptyset$ and a frame  $\{ \hat{e}^{\pm}_{i} \}$ on $V$ with $\nabla^{\pm} \hat{e}^{\pm}_{i} = 0$.
Then we have another Lagrange embedding
\[
g = g^{+}\oplus g^{-}:V \rightarrow \R^{n}_{\bby} \times \R^{n}_{\bbq}
\]
with different coordinates $\by = (y_{1},\cdots,y_{n})$, $\bq = (q_{1},\cdots,q_{n})$ satisfying
\[
dg = dg^{+}\oplus dg^{-} = \Phi^{+}\oplus\Phi^{-}
\]
and
\[
\frac{\rd}{\rd y_{i}} =\hat{e}_{i}^{+}, \quad \frac{\rd}{\rd q_{i}} =\hat{e}_{i}^{-}
\]
for $i\in \{ 1,\cdots,n \}$.
We show embeddings $f$ and $g$ are compatible with some affine Lagrange equivalences.
Since $\nabla^{\pm} e^{\pm}_{i} = \nabla^{\pm} \hat{e}^{\pm}_{i} = 0$, the base change $F$ from $\{ e^{\pm}_{i} \}$ to $\{ \hat{e}^{\pm}_{i} \}$ does not depend on points of $U\cap V$ and actually it is of the form
\[
F=
\begin{pmatrix}
  A & \bzero\\
  \bzero & (A^{\top})^{-1}
\end{pmatrix},
\]
where $A\in GL_{n}(\R)$.
It satisfies that $F \circ df = dg$ on $U\cap V$.
Set $\bb(p)=Af^{+}(p)-g^{+}(p)$ and $\bb^{\prime}(p) = (A^{\top})^{-1}f^{-}(p) - g^{-}(p)$ for $p\in U\cap V$.
Since $F\circ df =dg$, it holds that
\begin{align*}
  d\bb &= d(Af^{+} - g^{+}) = A\circ df^{+} - dg^{+} = \bzero,\\
  d\bb^{\prime} &= d((A^{\top})^{-1}f^{-} - g^{-}) = (A^{\top})^{-1} \circ df^{-} - dg^{-} = \bzero.
\end{align*}
Therefore $\bb$ and $\bb^{\prime}$ are constant functions.
We define an affine transformation $\Lcal_{1}:\R^{n}_{\bbx} \times \R^{n}_{\bbp}  \rightarrow \R^{n}_{\bby} \times \R^{n}_{\bbq}$ by
\[
\Lcal_{1}(\bx,\bp):= (A\bx+\bb,(A^{\top})^{-1}+\bb^{\prime}).
\]
Then it holds that $d\Lcal_{1} =F$.
Obviously $\Lcal_{1}\circ f = g$ holds on $U\cap V$ and $\Lcal_{1}$ is an affine Lagrange equivalence.
For each Lagrange embedding $f$ and $g$, we can lift it to a Legendre embedding of $\R^{2n+1}$ so that on $U\cap V$ the obtained Legendre submanifolds are mutually transformed by the affine Legendre equivalence $\Lcal$ by
\[
\Lcal(\bx,\bp,z):= (\Lcal_{1}(\bx,\bp),z+\bc^{\top} \bx + d)
\]
with $\bc = A^{\top} \bb^{\prime}$ and $d=0$.
By gluing such Legendre submanifolds, we obtain a quasi-Hessian structure on $M$.
\qed

\vskip\baselineskip

Finally we summarize our intrinsic characterization of quasi-Hessian manifold defined in \cite{Nakajima-Ohmoto}.

\begin{dfn}\upshape
  A quasi-Hessian structure is defined to be a quasi-Codazzi structure with dually flat connections.
\end{dfn}

\vskip\baselineskip
\section*{Acknowledgements}
The author would like to thank Prof. T. Ohmoto, his ex-supervisor, for guiding him to this subject and for instructions and discussions.
He is also grateful to Profs. H. Furuhata, A. Honda, G. Ishikawa and Dr. N. Nakajima for their useful comments at his talks in seminars.

\vskip\baselineskip
\vskip\baselineskip

\end{document}